\documentclass[11pt]{amsart}

\usepackage{amssymb,amsmath,color,hyperref}
\usepackage[mathscr]{eucal}

\theoremstyle{plain}
\newtheorem{thm}{Theorem}[section]
\newtheorem{theorem}[thm]{Theorem}

\theoremstyle{definition}

\newtheorem{ex}[thm]{Example}

\theoremstyle{remark}

\newtheorem{remark}{Remark}
\newtheorem*{remark*}{Remark}

\numberwithin{equation}{section}

        \newcommand{\field}[1]{{\mathbb{#1}}}
        \newcommand{\NN}{\field{N}}
        \newcommand{\ZZ}{\field{Z}}
        
        \newcommand{\RR}{\field{R}}

\begin{document}

\title[Monte Carlo methods on compact symplectic manifolds]{Monte Carlo methods on compact symplectic manifolds}

\author[Y. A. Kordyukov]{Yuri A. Kordyukov}
\address{Institute of Mathematics, Ufa Federal Research Centre, Russian Academy of Sciences, 112~Chernyshevsky str., 450008 Ufa, Russia} \email{yurikor@matem.anrb.ru}

%
%
\begin{abstract}
We build an unbiased Monte Carlo estimator of the integral of any $C^1$ function on a prequantized compact symplectic manifold against a smooth Riemannian volume form, taking for quadrature nodes the determinantal point process associated with an appropriate spectral projection of the Bochner-Schr\"odinger operator. We show that the estimator satisfies a central limit theorem, and the decay rate of the mean squared error reaches the optimal worst-case rate investigated by Bakhvalov in Euclidean spaces. These results extend previous results of Lemoine and Bardenet on Monte Carlo methods on compact complex manifolds. 
\end{abstract}


 \maketitle


\section{Introduction}
In \cite{LB24}, Lemoine and Bardenet propose a new randomized method for numerical integration
on a compact complex manifold with respect to a continuous volume form. Taking for quadrature nodes a suitable determinantal point process, they build an unbiased Monte Carlo estimator of the integral of any $C^1$ function, and show that the estimator satisfies a central limit theorem, with a faster rate than under independent sampling and previous DPP-based quadratures. In particular, the decay rate of the mean squared error for $N$ quadrature nodes reaches the optimal worst-case rate investigated by Bakhvalov \cite{B59} in Euclidean spaces. The determinantal point process used \cite{LB24} is characterized by the Bergman kernel of a holomorphic Hermitian line bundle. The considerations of this paper are heavily based on the work of Berman \cite{Berman14,Berman18}, especially on the central limit theorem \cite[Theorem 5.8]{Berman18}. We refer the reader to \cite{BH20,GBV19,HMBL23,LB24,PBD} and the references therein for more information on quadrature rules based on determinantal and repelled point processes, in particular, on manifolds.

The main goal of this note is to extend the results of \cite{LB24} to a more general class of compact manifolds, namely, to prequantized symplectic manifolds. For this purpose, we will use the determinantal point processes constructed in \cite{dpp}. They are associated with an appropriate spectral projection of the Bochner-Schr\"odinger operator. 

Let $X$ be a compact manifold and $dv_X$ a smooth volume form on $X$. We assume that $dv_X$ is the Riemannian volume form of a Riemannian metric $g$ on $X$. We assume that that $X$ is a symplectic manifold of even dimension $2n$ and denote by $\mathbf B$ the symplectic form on $X$. Furthermore, we assume that it satisfies the prequantization condition:
\[
[\mathbf B]\in H^2(X,2\pi\mathbb Z).
\]
Under these assumptions, we build an unbiased Monte-Carlo estimator of the integral of any $C^1$ function on $X$ against the volume form $dv_X$.

There exists a Hermitian line bundle $(L,h^L)$ on $X$ with a Hermitian connection $\nabla^L$ such that 
\begin{equation}\label{e:def-omega}
\mathbf B=iR^L, 
\end{equation} 
where $R^L$ is the curvature of the connection $\nabla^L $ defined as $R^L=(\nabla^L)^2$. 

For any $p\in \NN$, let $L^p:=L^{\otimes p}$ be the $p$th tensor power of $L$ and let
\[
\nabla^{L^p}: {C}^\infty(X,L^p)\to
{C}^\infty(X, T^*X \otimes L^p)
\] 
be the Hermitian connection on $L^p$ induced by $\nabla^{L}$. Consider the induced Bochner Laplacian $\Delta^{L^p}$ acting on $C^\infty(X,L^p)$ by
\begin{equation}\label{e:def-Bochner}
\Delta^{L^p}=\big(\nabla^{L^p}\big)^{\!*}\,
\nabla^{L^p},
\end{equation} 
where $\big(\nabla^{L^p}\big)^{\!*}: {C}^\infty(X,T^*X\otimes L^p)\to
{C}^\infty(X,L^p)$ is the formal adjoint of  $\nabla^{L^p}$ with respect to the $L^2$-inner products on $C^\infty(X,L^p)$ and ${C}^\infty(X,T^*X\otimes L^p)$ defined by the Riemannian volume form $dv_X$ and Hermitian structure $h^{L^p}$. 
Let $V\in C^\infty(X,\RR)$ be a real-valued function. 

Define the Bochner-Schr\"odinger operator $H_p$ acting on $C^\infty(X,L^p)$ by
\[
H_{p}=\frac 1p\Delta^{L^p}+V. 
\] 
Denote by $L^2(X,L^p)$ the Hilbert completion of $C^\infty(X,L^p)$ with respect to the $L^2$-inner product on $C^\infty(X,L^p)$ defined by the Riemannian volume form $dv_X$ and the induced Hermitian structure on $L^p$:
\begin{equation}\label{e:def-inner}
\langle s_1, s_2\rangle_{p}=\int_X\langle s_1(x),s_2(x)\rangle_{h^{L^p}} dv_X(x), \quad s_1, s_2\in C^\infty(X,L^p).
\end{equation}
The operator $H_p$ is essentially self-adjoint  in the Hilbert space $L^2(X,L^p)$ with initial domain  $C^\infty(X,L^p)$. 

\begin{remark}
Given a symplectic form $\mathbf B$ and a Riemannian metric $g$ on $X$, there exists a unique almost complex structure $J$ on the tangent bundle $TX$, which is compatible with $\mathbf B$ and $g$, meaning that $\mathbf B(Ju, Jv) = \mathbf B(u,v)$, $g(Ju, Jv) =g(u,v)$ for any $u,v\in TX$ and the associated symmetric bilinear form $g_{\mathbf B}(u, v) = \mathbf B(u,Jv)$ is a Riemannian metric. It is constructed as follows. 
For ${x}\in X$, let $B_{x} : T_{x}X\to T_{x}X$ be the skew-symmetric operator such that 
\begin{equation}
\label{e:Bx}
\mathbf B_{x}(u,v)=g(B_{x}u,v), \quad u,v\in T_{x}X. 
\end{equation}
Then $J$ is given by 
\[
J_{x}=B_{x}(-B_{x}^2)^{1/2}, \quad x\in X. 
\] 
The Riemannian metric $g_{\mathbf B}$ is called the Riemannian metric associated with $\mathbf B$ and $g$. It is clear that 
\[
g_{\mathbf B}(u,v)=g((B^*_xB_x)^{1/2}u,v), \quad u,v\in T_xX.
\]
If $g=g_{\mathbf B}$, then $(X,\mathbf B,g)$ is called an almost-K\"ahler manifold. 

In the setting of \cite{Berman14,Berman18,LB24}, $(X,\mathbf B,g)$ is assumed to be a K\"ahler manifold, which means that it is almost-K\"ahler and $J$ is integrable. In this case, by the Bochner-Kodaira-Nakano identity, we have 
\begin{equation}\label{e:BKN}
\Delta^{L^p}=2\Box^{L^p}+np,
\end{equation}
where $\Box^{L^p}=\bar\partial^{L^p*}\bar\partial^{L^p}$ is the Kodaira-Laplacian on functions.
\end{remark}

\begin{remark}
 Assume that the Hermitian line bundle $(L,h^L)$ is trivial. Then we can write $\nabla^L=d-i \mathbf A$ with a real-valued 1-form $\mathbf A$ (the magnetic potential), and we have
\[
R^L=-id\mathbf A,\quad \mathbf B=d\mathbf A. 
\]
The operator $H_p$ is related with the semiclassical magnetic Schr\"odinger operator
\[
H_p=\hbar^{-1}[(i\hbar d+\mathbf A)^*(i\hbar d+\mathbf A)+\hbar V], \quad \hbar=\frac{1}{p},\quad p\in \NN.
\]
So $\mathbf A$ can be interpreted as the magnetic potential, $\mathbf B$ as the magnetic field and $V$ as the electric potential. 
\end{remark} 


As shown in \cite{Kor22} (see also \cite{charles20b}), the spectrum of $H_{p}$ as $p\to \infty$ admits the following asymptotic description.

For ${x}\in X$, denote by $\pm i a_j(x), j=1,\ldots,n,$ with $a_j(x)>0$ the eigenvalues of the skew-symmetric operator $B_{x} : T_{x}X\to T_{x}X$ defined by \eqref{e:Bx}. 
For any $\mathbf k=(k_1,\cdots,k_n)\in\ZZ_+^n$ and $x\in X$, set
\begin{equation}\label{e:def-Lambda}
\Lambda_{\mathbf k}(x)=\sum_{j=1}^n(2k_j+1) a_j(x)+V(x).
\end{equation}
These numbers can be interpreted as eigenvalues of infinite multiplicity of the magnetic Schr\"odinger operator in $T_{x}X$ with constant magnetic field and electric potential and can be naturally called the local Landau levels at $x$. 

Set
\begin{equation}\label{e:def-Sigma}
\Sigma=\left\{\Lambda_\mathbf {k}(x)\,:\, \mathbf k\in\ZZ_+^n, x\in X \right\}.
\end{equation}

\begin{theorem}[\cite{Kor22}, Theorem 1]\label{t:spectrum}
For any $K>0$, there exists $c>0$ such that for any $p\in \NN$ the spectrum of $H_{p}$ in the interval  $[0,K]$  is  contained in the $cp^{-1/4}$-neighborhood of $\Sigma$.  
\end{theorem}

The set $\Sigma$ is a closed subset of $\RR$, which can be represented as the union of closed intervals (bands):
\[
\Sigma=\bigcup_{\mathbf k\in\ZZ_+^n}[\alpha_{\mathbf k}, \beta_{\mathbf k}]
\]
where, for any $\mathbf k\in\ZZ_+^n$, the interval $[\alpha_{\mathbf k}, \beta_{\mathbf k}]$ is the image of the continuous function $\Lambda_{\mathbf k}$ on $X$: $[\alpha_{\mathbf k}, \beta_{\mathbf k}]=\{\Lambda_{\mathbf k}({x}) : x\in X\}$.
In general, the bands $[\alpha_{\mathbf k}, \beta_{\mathbf k}]$ may overlap without any gaps so that $\Sigma$ is the semi-axis $[\Lambda_0,+\infty)$ with $\Lambda_0=\inf_{x\in X} \Lambda_0(x)$, where $\Lambda_0(x)$ is the lowest Landau level at $x$: 
\[
\Lambda_0(x):=\sum_{j=1}^n a_j(x)+ V(x). 
\]

There are some cases when $\Sigma$ has a gap, that is, $[\Lambda_0,+\infty)\setminus \Sigma \neq \emptyset$. Then Theorem~\ref{t:spectrum} predicts the existence of gaps in the spectrum of the operator $H_{p}$ for any $p$ large enough. 

\begin{ex}
Note that $(X,\mathbf B,g)$ is an almost-K\"ahler manifold iff $B=J$ or equivalently $a_j(x)\equiv 1$, $j=1,2,\ldots,n$. In this case, if $V(x)\equiv 0$, then 
\[
\Lambda_{\mathbf k}(x)=2|\mathbf k|+n, \quad \mathbf k\in\ZZ_+^n,  
\]
and $\Sigma=\{2N+n:N\in \mathbb Z_+\}$ is a countable discrete set. In particular, this shows that, for a  given $(X,\mathbf B)$, we can always choose $g$ and $V$ such that $\Sigma$ has a gap. 

More generally, we may assume that the functions $a_j$ are constants. The set $\Sigma$ will also have gaps if $a_j$'s are not constants, but vary slowly enough.
\end{ex}

\begin{ex}
For an arbitrary $(X,\mathbf B,g)$, one can always construct the operator $H_p$ with gaps in $\Sigma$, choosing 
\[
V(x)=-\sum_{j=1}^n a_j(x)=\frac{1}{2}\operatorname{tr}((-B_{x}^2)^{1/2}), \quad x\in X.
\]
From the last equality and uniform positivity of $B$, it follows that $V$ is a smooth function. The corresponding operator 
\[
\Delta_p:=pH_p=\Delta^{L^p}+pV
\] 
was introduced by Guillemin and Uribe in \cite{Gu-Uribe}. It is called the renormalized Bochner Laplacian. 
For this operator, $\Lambda_0=0$ is an isolated point in $\Sigma$. 

When $(X,\mathbf B,g)$ is a K\"ahler manifold, by \eqref{e:BKN}, the renormalized Bochner Laplacian $\Delta_p$ is twice the Kodaira-Laplacian on functions $\Box^{L^p}=\bar\partial^{L^p*}\bar\partial^{L^p}$. In this case, $\Lambda_0=0$ is an eigenvalue of $\Delta_p$ for any $p$, and the corresponding eigenspace consists of holomorphic sections of $L^p$. 

In the general case of a symplectic manifold, Guillemin and Uribe in \cite{Gu-Uribe} suggest to consider the spectral subspace of the renormalized Bochner Laplacian $\Delta_p$ corresponding to sufficiently small neighborhood of zero as a substitute of the space of holomorphic sections of $L^p$.
\end{ex}

We will assume that $\Sigma$ has a gap and take an interval $I=(\alpha,\beta)$ such that $\alpha,\beta\not \in \Sigma$. 
By Theorem \ref{t:spectrum}, there exists $\mu_0>0$ and $p_0\in \NN$ such that for any $p>p_0$ 
\[
\sigma(H_{p})\subset (-\infty, \alpha-\mu_0) \cup I \cup (\beta+\mu_0, \infty).
\] 
Let $P_{p,I}$ be the spectral projection of $H_{p}$ associated with $I$:
\[
P_{p,I} : L^2(X,L^p)\to \mathcal H_p=\operatorname{Im} P_{p,I}. 
\]
We will consider the manifold $X$ equipped with the Riemannian volume form $dv_X$ and the determinantal point process on $X$ associated with the finite rank projection $P_{p,I}$ for $p>p_0$. It can be constructed as follows. 

Set $N_p:=\dim \mathcal H_p$.  Consider the $N_p$-fold product $X^{N_p}$, which can be viewed as the configuration space of $N_p$ particles on $X$. 
Let $\{s_{p,j} \in \mathcal H_p: j=1,\ldots, N_p\}$ be an orthonormal basis in $\mathcal H_p$. We define the Slater determinant as the section $\Psi_p \in C^\infty(X^{N_p}, (L^p)^{\boxtimes N_p})$ of the line bundle $(L^p)^{\boxtimes N_p}$ over $X^{N_p}$ given for any $(x_1, x_2, \ldots, x_{N_p}) \in X^{N_p}$ by
\begin{equation}\label{e: Slater}
\Psi_p(x_1, x_2, \ldots, x_{N_p}) := \det(s_{p,j}(x_i))_{i,j=1}^{N_p}, 
\end{equation}
which does not depend on the choice of orthonormal basis of $\mathcal H_p$.

The determinantal point process on $X$ associated with $P_{p,I}$ is defined by the probability measure $d\nu_{N_p}$ on $X^{N_p}$ given by
\[
d\nu_{N_p}:=\frac{1}{N_p!} |\Psi_p|^2_{h^{L^p}} dv_X^{N_p}.  
\]
In \cite{dpp}, we studied the asymptotic behavior of smooth linear statistics for this process as $p$ goes to infinity (In the current setting this limit corresponds to the limit of a large number $N_p$ of particles). As a consequence, we obtained the law of large numbers and central limit theorem. 
 



Before we state the main result of the paper, let us introduce some notation. 

Let 
$$
\mathcal K_I:=\{\mathbf k\in \ZZ^n_+ : \Lambda_{\mathbf k}(x)\in I\}.
$$
Under current assumption on $I$, this set is independent of $x\in X$.

Fix $x\in X$.  Choose an orthonormal base $\{e_j : j=1,\ldots,2n\}$ in $T_{x}X$ such that  
\[
B_{x}e_{2k-1}=a_k(x)e_{2k}, \quad B_{x}e_{2k}=-a_k(x)e_{2k-1},\quad k=1,\ldots,n. 
\]

For any $m=1,\ldots,n$, introduce a function $I_m : \mathbb Z^n_+\times \mathbb Z^n_+\to \mathbb Z$, setting, for $\mathbf k,\mathbf k^\prime\in \mathbb Z^n_+$,
\begin{equation}
\begin{aligned} \label{e:defIm}
I_m(\mathbf k,\mathbf k)& =2k_m+1,\\
I_m(\mathbf k+\epsilon_j,\mathbf k)=I_m(\mathbf k,\mathbf k+\epsilon_j)& =-(k_j+1)\delta_{jm}, \quad j=1,\ldots,n, \\
I_m(\mathbf k,\mathbf k^\prime)& =0, \quad |\mathbf k-\mathbf k^\prime|>1,
\end{aligned}
\end{equation}
where $(\epsilon_1,\ldots,\epsilon_n)$ is the standard basis in $\mathbb Z^n$. 

Set                                             
\[
\alpha_m=\sum_{\mathbf k^\prime,\mathbf k^{\prime\prime}\in \mathcal K_I} I_m(\mathbf k^\prime,\mathbf k^{\prime\prime}).
\]                  
             
For any $f \in C^1(X)$, define  
\begin{equation}\label{e:def-dfI}
|df(x)|_{I}^2= \sum_{m=1}^{n} \frac{\alpha_m}{a_m(x)}   \left[\left(\nabla_{e_{2m-1}}f(x)\right)^2+\left(\nabla_{e_{2m}}f(x)\right)^2\right],
\end{equation}
where, for any $v\in T_xX$, $\nabla_vf(x)$ denotes the derivative of $f$ in the direction of $v$. 

The function $x\in X\mapsto |df(x)|_{I}^2$ is a well-defined continuous function on $X$. This follows from an interpretation of this function in terms of the spectral data of the model operators associated with the Bochner-Schr\"odinger operator $H_p$ given in \cite{dpp} (see Section 4, in particular, formula (4.3)).  

\begin{remark}
In the case when $|\mathcal K_I|=1$, that is, $\mathcal K_I$ consists of a single element $\mathbf k\in \mathbb Z_+^n$, we have 
\[
\alpha_m=I_m(\mathbf k,\mathbf k)=2k_m+1. 
\]
It follows that
\[
|df(x)|_{I}^2= \sum_{m=1}^{n} \frac{2k_m+1}{a_m(x)}   \left[\left(\nabla_{e_{2m-1}}f(x)\right)^2+\left(\nabla_{e_{2m}}f(x)\right)^2\right].
\]

When $\mathcal K_I$ consists of a single element $\mathbf k=0$, the function has the following geometric interpretation \cite{dpp}:
\begin{equation}\label{e:dfI-ex}
|df(x)|_{I}^2=|df(x)|_{g^{-1}_{\mathbf B}}^2,
\end{equation}
where $g_{\mathbf B}$ is the Riemannian metric associated with $\mathbf B$ and $g$ and $g^{-1}_{\mathbf B}$ is the induced Riemannian metric on $T^*X$. 
\end{remark}

%
%
%
%

Recall that $P_{p,I}$ denotes the spectral projection of $H_{p}$ associated with the interval $I$. Let $P_{p,I}(x,x^\prime)$, $x,x^\prime\in X$, be its smooth Schwartz kernel with respect to the Riemannian volume form $dv_X$. 

Let $\Omega_{\mathbf B}=\frac{1}{n!} \mathbf B^n$ be the Liouville volume form on $X$. 
For any function $f$ on $X$, define a function $f_{\mathbf B}$ on $X$ by $$f_{\mathbf B}(x)=\frac{f(x)}{\sqrt{\det B_x}}, \quad x\in X.$$ Thus, $fdv_X=f_{\mathbf B}\Omega_{\mathbf B}$. 

The main result of the paper is the following theorem. 

\begin{thm}\label{t:mcarlo}  
For any $f \in C^1(X, \mathbb R)$, the random variable $\Xi_p$ on $(X^{N_p}, d\nu_{N_p})$ given by 
\begin{multline}\label{e:Xip}
\Xi_p(x_1, x_2, \ldots, x_{N_p})=N^{\frac{1}{2n}+\frac 12}_p\left(\sum_{i=1}^{N_p}\frac{f(x_i)}{P_{p,I}(x_i,x_i)}-\int_X f(x)dv_X(x) \right),\\ (x_1, x_2, \ldots, x_{N_p})\in X^{N_p}.
\end{multline}
converges in distribution as $p\to \infty$ to a centered normal random variable $N(0,\sigma^2)$ with variance 
\begin{equation}\label{e:mcarlo}
\sigma^2=\frac{1}{4\pi}\frac{({\rm Vol}_{\mathbf B}(X))^{(n+1)/n}}{|\mathcal K_I|^{(n-1)/n}}\int_X|df_{\mathbf B}(x)|_{I}^2\,\Omega_{\mathbf B}(x).
\end{equation}
\end{thm}

This theorem allows us to build a Monte Carlo estimator for the integral $J=\int_X fdv_X$ by the formula
\[
\hat{J}(x_1,\ldots,x_{N_p})=\sum_{i=1}^{N_p}\frac{f(x_i)}{P_{p,I}(x_i,x_i)},\quad (x_1, x_2, \ldots, x_{N_p})\in X^{N_p}. 
\]
It is unbiased, because by a well-known formula for the expectation of linear statistics for determinantal point processes (cf. \eqref{e:EpXip} below), we immediately get 
\[
\mathbb E[\hat{J}]=\int_X f(x) dv_X(x).
\]
Therefore, the mean square error of the estimator equals its variance, which is asymptotically of order $N_p^{-1-2/2n}$:
\[
\operatorname{MSE}\sim \frac{\sigma^2}{N^{\frac{n+1}{n}}_p}\quad p\to\infty. 
\]
This matches the optimal worst-case rate established by N. S. Bakhvalov \cite{B59} for randomized integration of $C^1$ functions in Euclidean spaces; see also \cite{N16}. 

The important difference of our setting from the one of Lemoine and Bardenet \cite{LB24} and Berman \cite{Berman14,Berman18} is that, in our case, the space $\mathcal H_p$ depends on the metric $g$ (in \cite{LB24}, it is the space of holomorphic sections of $L^p$). By this reason, the determinantal point process depends on the metric as well. Of course, one can choose $g=g_{\mathbf B}$, then this reduces to the choice of a positive almost complex structure compatible with $\mathbf B$.  
On the other hand, even in the K\"ahler case as in \cite{LB24} our setting is broader, because the holomorphic sections are associated with the lowest Landau levels of the magnetic Laplacian, whereas we can consider the eigenspaces associated with higher Landau levels (the so-called polyanalytic functions).


Let us consider an example discussed in \cite{dpp} (see \cite[Section 5]{dpp} for more details). Assume that $g=g_{\mathbf B}$, or equivalently $a_1(x)=a_2(x)=\ldots=a_n(x)=1$ for any $x\in X$ and $V(x)\equiv 0.$
Then 
\[
\Sigma=\{2N+n:N\in \mathbb Z_+\}.
\]
Take an interval around $N$th Landau level, say, $$I_N=(2N+n-1,2N+n+1).$$ 
For the Landau Hamiltonian, this case corresponds to the pure $N+1$-analytic Ginibre point process (see \cite{HW19}). 

We compute
\[
\mathcal K_{I_N}:=\{\mathbf k\in \ZZ^n_+ : |\mathbf k|=N\}, \quad
|\mathcal K_{I_N}|=\binom{N+n-1}{n-1}=\frac{(N+n-1)!}{N!(n-1)!},
\]
and
\[
\alpha_m=\frac{2N+n}{n}|\mathcal K_{I_N}|, \quad m=1,\ldots,n.
\]
By \eqref{e:def-dfI}, we have
\[
|df(x)|_{I_N}^2=\frac{2N+n}{n}|\mathcal K_{I_N}||df(x)|^2_{g^{-1}_{\mathbf B}},
\]
and, by \eqref{e:mcarlo},
\[
\sigma^2=\frac{1}{4\pi}\frac{2N+n}{n}|\mathcal K_{I_N}|^{1/n}({\rm Vol}_{\mathbf B}(X))^{(n+1)/n}\int_X|df_{\mathbf B}(x)|^2_{g^{-1}_{\mathbf B}} \,\Omega_{\mathbf B}(x).
\]


Take $I=(n-1,2N+n+1)$. For the Landau Hamiltonian, this case corresponds to the full $N+1$-analytic Ginibre point process (see \cite{HW19}). We compute
\[
\alpha_m=\sum_{j=0}^N|\mathcal K_{I_j}|=|\mathcal K_{I}|=\binom{N+n}{n}=\frac{N+n}{n}|\mathcal K_{I_N}|.
\]
By \eqref{e:def-dfI}, we have
\[
|df(x)|_{I}^2=\frac{N+n}{n}|\mathcal K_{I_N}||df(x)|^2_{g^{-1}_{\mathbf B}}.
\]
and, by \eqref{e:mcarlo}, 
\[
\sigma^2=\frac{1}{4\pi}\left(\frac{N+n}{n}\right)^{1/n}|\mathcal K_{I_N}|^{1/n} ({\rm Vol}_{\mathbf B}(X))^{(n+1)/n}\int_X |df_{\mathbf B}(x)|^2_{g^{-1}_{\mathbf B}}\,\Omega_{\mathbf B}(x).
\]
Thus, we see the estimate for the mean square error is better in the second case, which is due to the repelling property of determinantal point processes.  

In the general case, we have the following upper estimate for $\sigma^2$: 
\begin{equation}\label{e:upper-est}
\sigma^2\leq C|\mathcal K_I|^{1/n} ({\rm Vol}_{\mathbf B}(X))^{(n+1)/n}  \int_X|df_{\mathbf B}(x)|_{g^{-1}_{\mathbf B}}^2\,\Omega_{\mathbf B}(x),
\end{equation}
with some constant $C>0$, independent of $f$. 

Indeed, denoting $I=(\alpha,\beta)$ and 
$$
a_{min}=\min \{a_j(x) : x\in X, j=1,2,\ldots,n\}, \quad V_{min}= \min \{V(x):x\in X\}, 
$$
we get, for any $\mathbf k\in \mathcal K_I$ and $m=1,2,\ldots,n$, 
\[
(2k_m+1)a_{min}+V_{min}\leq \Lambda_{\mathbf k}(x)=\sum_{j=1}^n(2k_j+1) a_j(x)+V(x)\leq \beta, 
\]
and, as a consequence, 
\[
 I_m(\mathbf k,\mathbf k)=2k_m+1\leq \frac{\beta-V_{min}}{a_{min}}. 
\]
Using that $I_m(\mathbf k^\prime,\mathbf k^{\prime\prime})\leq 0$ for $\mathbf k^\prime\neq \mathbf k^{\prime\prime}$, we infer that 
\[
\alpha_m=\sum_{\mathbf k^\prime,\mathbf k^{\prime\prime}\in \mathcal K_I} I_m(\mathbf k^\prime,\mathbf k^{\prime\prime})\leq \sum_{\mathbf k\in \mathcal K_I} I_m(\mathbf k,\mathbf k)\leq \frac{\beta-V_{min}}{a_{min}}|\mathcal K_I|.
\]  
By \eqref{e:def-dfI} and \eqref{e:dfI-ex}, we have
\begin{align*}
|df(x)|_{I}^2\leq & \frac{\beta-V_{min}}{a_{min}} |\mathcal K_I| \sum_{m=1}^{n} \frac{1}{a_m(x)}   \left[\left(\nabla_{e_{2m-1}}f(x)\right)^2+\left(\nabla_{e_{2m}}f(x)\right)^2\right]\\
\leq & \frac{\beta-V_{min}}{a_{min}}|\mathcal K_I| |df(x)|_{g^{-1}_{\mathbf B}}^2.
\end{align*}
and, by \eqref{e:mcarlo},
\[
\sigma^2\leq \frac{1}{4\pi}\frac{\beta-\min V(x)}{\min a_1(x)}|\mathcal K_I|^{1/n} ({\rm Vol}_{\mathbf B}(X))^{(n+1)/n}  \int_X|df_{\mathbf B}(x)|_{g^{-1}_{\mathbf B}}^2\,\Omega_{\mathbf B}(x),
\]
that proves \eqref{e:upper-est}.

In the next section, we give the proof of the main result, Theorem~\ref{t:mcarlo}. 
 
\section{Proof of Theorem~\ref{t:mcarlo}}

The proof of Theorem~\ref{t:mcarlo} closely follows the proof of Theorem 1.2 in \cite{LB24}, which is in its turn based on the proof of the central limit theorem for linear statistics proved by Berman, see \cite[Theorem 1.5]{Berman18}. There are two essential differences. We use more refined estimates for the kernel of the spectral projection obtained by techniques developed by Ma and Marinescu \cite{ma-ma:book,ma-ma08}  and extended to the current setting by the author in \cite{Kor18,Kor22}. In \cite{dpp}, this allows us to give a simpler proof of an analog of \cite[Theorem 5.8]{Berman18} in the current setting. Second, we use a second order differential operator to define  the space $\mathcal H_p$ and therefore the space $\mathcal H_p$ depends on the Hermitian structure of the line bundle $L^p$ unlike the Bergman space used in \cite{Berman18,LB24}. Therefore, the trick by Berman based on the use an auxiliary line bundle doesn't work, and we apply some operator theoretic methods instead.   

Given $f \in C^1(X, \mathbb R)$, for any $p\in \mathbb N$, we introduce functions $u_p$ and $f_p$ on $X$ by
 \begin{equation}\label{e:def-fp}
f_p(x)=\frac{N_p}{P_{p,I}(x,x)}f(x), \quad x\in X, 
 \end{equation}
and 
 \begin{equation}\label{e:def-up}
u_p(x)=N^{\frac{1}{2n}-\frac 12}_p\left(f_p(x)-\int_X fdv_X \right), \quad x\in X.
 \end{equation}


Recall that, for any function $u : X\to \mathbb R$ and for any $p \in \mathbb N$, the associated linear statistics is the random variable over $(X^{N_p}, d\nu_{N_p})$ defined by 
\[
\mathcal N_p[u](x_1, x_2, \ldots, x_{N_p})=\sum_{j=1}^{N_p}u(x_j),\quad (x_1, x_2, \ldots, x_{N_p})\in X^{N_p}.
\]
It is easy to see that the random variable $\Xi_p$ defined by \eqref{e:Xip} is the linear statistic associated with $u_p$: 
\[
\Xi_p=\mathcal N_p[u_p].
\]

Consider a family of smooth volume forms on $X$ given by $$\mu_t=e^{-tu_p} dv_X, \quad t\in \mathbb R.$$
Denote by $\langle \cdot, \cdot\rangle_{p,t}$, $t\in \mathbb R$, the $L^2$-inner product on $L^2(X,L^p)$ defined by the smooth volume form $\mu_t$ and Hermitian structure on $L^p$:
\[
\langle s_1, s_2\rangle_{p,t}=\int_X\langle s_1(x),s_2(x)\rangle_{h^{L^p}} e^{-tu_p(x)} dv_X(x), \quad s_1, s_2\in L^2(X,L^p).
\]
For $t=0$, it coincides with the inner product $\langle \cdot, \cdot\rangle_{p}$ defined by \eqref{e:def-inner}. We will denote by $L^2(X,L^p;\mu_t)$ the space $L^2(X,L^p)$ equipped with the inner product $\langle \cdot, \cdot\rangle_{p,t}$. 

Let $P_{p,I,t}$ be the orthogonal projection in $L^2(X,L^p;\mu_t)$ on $\mathcal H_p$ and $P_{p,I,t}\in C^\infty(X\times X, L^p\boxtimes (L^p)^*)$ its integral kernel with respect to $\mu_t$:
\[
P_{p,I,t}s(x)=\int_X P_{p,I,t}(x,y)s(y)e^{-tu_p(y)} dv_X(y), \quad s\in L^2(X,L^p;\mu_t). 
\]

Let $(X^{N_p}, d\nu_{p,t})$ be the determinantal point process associated with $P_{p,I,t}$. 
One can show that
\begin{equation}\label{e:dnupt}
d\nu_{p,t}=\frac{|\Psi_{p}|^2_{h^{L^p}} e^{-t\Xi_p} dv_X^{N_p}}{\int_{X^{N_p}} |\Psi_{p}|^2_{h^{L^p}} e^{-t\Xi_p} dv_X^{N_p}},
\end{equation}
where $\Psi_{p}$ is given by \eqref{e: Slater}.

Indeed, let $\{s_{j,p,t} \in \mathcal H_p: j=1,\ldots, N_p\}$ be an orthonormal basis in $\mathcal H_p$ with respect to $\langle\cdot,\cdot\rangle_{p,t}$ and $\Psi_{p,t} \in C^\infty(X^{N_p}, (L^p)^{\boxtimes N_p})$ the corresponding Slater determinant:
\[
\Psi_{p,t}(x_1, x_2, \ldots, x_{N_p}) := \det(s_{j,p,t}(x_i))_{i,j=1}^{N_p}. 
\]
By definition, the measure $d\nu_{p,t}$ on $X^{N_p}$ given by
\begin{equation}\label{e:dnupt1}
d\nu_{p,t}:=\frac{1}{N_p!} |\Psi_{p,t}|^2_{h^{L^p}} e^{-t\Xi_p} dv_X^{N_p}. 
\end{equation}
We can write
\[
s_{j,p,t}=\sum_{k=1}^{N_p}\langle s_{j,p,t}, s_{p,k}\rangle_{p}s_{p,k}, \quad  j=1,\ldots, N_p.
\]
It follows that
\begin{equation}\label{e:dnupt3}
\Psi_{p,t}=(\det T)\Psi_{p},
\end{equation}
where $T=(\langle s_{j,p,t}, s_{p,k}\rangle_{p})_{j,k=1}^{N_p}\in {\rm GL}(N_p,\mathbb C)$,
and, by \eqref{e:dnupt1},
\begin{equation}\label{e:dnupt4}
d\nu_{p,t}=\frac{|\det T|^2}{N_p!} |\Psi_{p}|^2_{h^{L^p}} e^{-t\Xi_p} dv_X^{N_p}. 
\end{equation}
In particular, we have
\begin{equation}\label{e:dnupt5}
\int_{X^{N_p}}d\nu_{p,t}=\frac{|\det T|^2}{N_p!} \int_{X^{N_p}} |\Psi_{p}|^2_{h^{L^p}} e^{-t\Xi_p} dv_X^{N_p}=1. 
\end{equation}
From \eqref{e:dnupt4} and \eqref{e:dnupt5}, we immediately get \eqref{e:dnupt}. 
 

Denote by $\mathbb E_{p,t}$ the expectation defined by $(X^{N_p}, d\nu_{p,t})$:
\[
\mathbb E_{p,t}[F]:=\int_{X^{N_p}}F d\nu_{p,t},\quad F\in L^1(X^{N_p}, d\nu_{p,t}).
\]
In particular,
\[
\mathbb E_{p,0}[F]=\mathbb E_{p}[F] :=\int_{X^{N_p}}F d\nu_{N_p},\quad F\in L^1(X^{N_p}, d\nu_{N_p}).
\] 
By \eqref{e:dnupt1}, we have
\begin{equation}\label{e:Eptf}
\mathbb E_{p,t}[F]=\frac{\int_{X^{N_p}}F |\Psi_{p}|^2_{h^{L^p}} e^{-t\Xi_p} dv_X^{N_p}}{\int_{X^{N_p}} |\Psi_{p}|^2_{h^{L^p}} e^{-t\Xi_p} dv_X^{N_p}}=\frac{\mathbb E_p[Fe^{-t\Xi_p}]}{\mathbb E_p[e^{-t\Xi_p}]}.
\end{equation}

We consider the log Laplace transform 
\[
F_p(t)=-\log \mathbb E_p[e^{-t\Xi_p}], \quad t\in\mathbb R.
\]
It is clear that
\begin{equation}\label{e:Fp0}
F_p(0)=0.
\end{equation}
Using \eqref{e:Eptf}, we compute
\[
\frac{dF_p}{dt}(t)=\frac{\mathbb E_p[\Xi_pe^{-t\Xi_p}]}{\mathbb E_p[e^{-t\Xi_p}]}=\mathbb E_{p,t}[\Xi_p], \quad t\in\mathbb R.
\]
We observe that
\begin{equation}\label{e:dFp0}
\frac{dF_p}{dt}(0)=0.
\end{equation} 
Indeed, using \eqref{e:def-up} and a well-known formula for the expectation of linear statistics for determinantal point processes, we have
\begin{multline}\label{e:EpXip}
\mathbb E_p[\Xi_p]=\int_XP_{p,I}(x,x)u_p(x)dv_X(x)\\
=N^{\frac{1}{2n}-\frac 12}_p(N_p\int_Xf(x)dv_X(x)-\int_XP_{p,I}(x,x)dv_X(x)\int_Xf(x)dv_X(x))=0.
\end{multline}


Next, we compute the second derivative of $F_p$:
\[
\frac{d^2F_p}{dt^2}(t)=-\frac{\mathbb E_p[\Xi^2_pe^{-t\Xi_p}]}{\mathbb E_p[e^{-t\Xi_p}]}+\frac{\mathbb E_p[\Xi_pe^{-t\Xi_p}]^2}{\mathbb E_p[e^{-t\Xi_p}]^2}
=-\mathbb E_{p,t}[\Xi^2_p]+\mathbb E_{p,t}[\Xi_p]^2=-\mathbb V_{p,t}[\Xi_p],
\]  
where $\mathbb V_{p,t}[\Xi_p]$ is the variance of the random variable $\Xi_p$ on $(X^{N_p}, d\nu_{p,t})$. 

By a well-known formula for the variance of linear statistics, we have
\[
\frac{d^2F_p}{dt^2}(t)=-\frac{1}{2}\int_X\int_X |P_{p,I,t}(x,y)|_{h^{L^p}}^2(u_p(x)-u_p(y))^2d\mu_t(x)d\mu_t(y).
\]
Using \eqref{e:def-up}, we can rewrite this formula as follows:
\begin{multline}\label{e:d2F1}
\frac{d^2F_p}{dt^2}(t)\\ =-N^{\frac{1}{n}-1}_p\frac{1}{2}\int_X\int_X |P_{p,I,t}(x,y)|_{h^{L^p}}^2(f_p(x)-f_p(y))^2d\mu_t(x)d\mu_t(y).
\end{multline}

For an integral operator $K$ in $L^2(X,L^p)$ with smooth integral kernel, denote by $\|K\|_{t,2}$ its Hilbert-Schmidt norm as an operator in the Hilbert space $L^2(X,L^p;\mu_t)$. We have 
\[
\|K\|_{t,2}^2=\int_X \int_X |K^{(t)}(x, y)|^2_{h^{L^p}}d\mu_t(x)d\mu_t(y),
\]
where $K^{(t)}\in C^\infty(X\times X, L^p\boxtimes (L^p)^*)$ is the integral kernel of $K$ with respect to $d\mu_t$:
\[
Ks(x)=\int_X K^{(t)}(x,y)s(y)d\mu_t(y),\quad x\in X, s\in C^\infty(X,L^p).
\]
It is clear that 
\[
K^{(t)}(x,y)=K(x,y)e^{tu_p(y)}, \quad x,y\in X,
\]
where $K\in C^\infty(X\times X, L^p\boxtimes (L^p)^*)$ is the integral kernel of $K$ with respect to $dv_X$.
It follows that 
\begin{equation}\label{e:HS-Kt}
\|K\|_{t,2}=\|e^{-(t/2)u_p}Ke^{(t/2)u_p}\|_{2}. 
\end{equation}

It is easy to see that the integral in the left-hand side of \eqref{e:d2F1} is written as 
\begin{multline}\label{e:VNpft-rhs}
\int_X \int_X |P_{p,I,t}(x, y)|^2_{h^{L^p}} (f_p(x)-f_p(y))^2d\mu_t(x)d\mu_t(y)\\ =\|[P_{p,I,t},f_p]\|^2_{t,2}
\end{multline} 

Now we consider the limit $p\to \infty$.  We know that 
\[
|e^{-tu_p(x)}-1|\leq C_1t\sup_{x\in X} |u_p(x)|\leq C_2tN^{\frac{1}{2n}-\frac 12}_p\leq C_3tp^{-\frac{n-1}{2}}, \quad x\in X.
\]
It follows that 
\begin{equation}\label{e:exp}
\|e^{-tu_p}-1\|\leq Ctp^{-\frac{n-1}{2}},
\end{equation}
where, for a bounded operator $A$ in $L^2(X,L^p)$, we denote by $\|A\|$ its operator norm. 

By \eqref{e:HS-Kt} and \eqref{e:exp}, we have
\begin{equation}\label{e:com-Ppit22}
\|[P_{p,I,t},f_p]\|_{t,2}=\|e^{-(t/2)u_p}[P_{p,I,t},p]e^{(t/2)u_p}\|_{2}=\|[P_{p,I,t},f_p]\|_{2}(1+\mathcal O(tp^{-\frac{n-1}{2}})).
\end{equation}

Now we claim that
\begin{equation}\label{e:com-Ppit}
\|[P_{p,I,t},f_p]\|_{2}=\|[P_{p,I},f_p]\|_{2}(1+\mathcal O(tp^{-\frac{n-1}{2}})), \quad p\to\infty.  
\end{equation}

If $p$ is large enough, the operator ${\rm Id}+P_{p,I}(e^{-tu_p}-1)P_{p,I}$ is invertible as a bounded operator in $L^2(X,L^p)$ and
\begin{equation}\label{e:Ppit}
P_{p,I,t} = ({\rm Id}+P_{p,I}(e^{-tu_p}-1)P_{p,I})^{-1} P_{p,I} e^{-tu_p}.
\end{equation}
Indeed, for any $s\in L^2(X,L^p)$, we have 
\[
\langle e^{-tu_p}(s - P_{p,I,t}s), s_1\rangle_{p}=\langle s - P_{p,I,t} s, s_1\rangle_{p,t}= 0, \quad s_1 \in \mathcal H_p.
\]
It follows that 
\[
P_{p,I} e^{-tu_p}(s - P_{p,I,t})s=0
\]
or
\[
P_{p,I} e^{-tu_p}P_{p,I,t}=P_{p,I} e^{-tu_p}. 
\]
Since $P_{p,I,t}=P_{p,I} P_{p,I,t}$, we can write 
\[
(P_{p,I} e^{-tu_p}P_{p,I}+{\rm Id}-P_{p,I})P_{p,I,t}=P_{p,I} e^{-tu_p},
\]

By \eqref{e:exp}, if $p$ is large enough, the operator ${\rm Id}+P_{p,I}(e^{-tu_p}-1)P_{p,I}$ is invertible in $L^2(X,L^p)$ and
\begin{equation}\label{e:inverse}
\|({\rm Id}+P_{p,I}(e^{-tu_p}-1)P_{p,I})^{-1}\|=1+\mathcal O(tp^{-\frac{n-1}{2}}), \quad p\to \infty. 
\end{equation}
In particular, this gives the desired identity \eqref{e:Ppit}. 

By \eqref{e:Ppit}, we get
\[
[P_{p,I,t},f_p] = ({\rm Id}+P_{p,I}(e^{-tu_p}-1)P_{p,I})^{-1} [P_{p,I},f_p] e^{-tu_p}
\]
\[
-({\rm Id}+P_{p,I}(e^{-tu_p}-1)P_{p,I})^{-1}([P_{p,I},f_p](e^{-tu_p}-1)P_{p,I}+P_{p,I}(e^{-tu_p}-1) 
[P_{p,I},f_p])
\]
\[
\times ({\rm Id}+P_{p,I}(e^{-tu_p}-1)P_{p,I})^{-1} P_{p,I} e^{-tu_p}.
\]
Using this identity and the estimates \eqref{e:exp} and \eqref{e:inverse}, we can easily complete the proof of \eqref{e:com-Ppit}.

By \cite[Theorem 3]{Kor22}, we know that 
\begin{equation}\label{e:PpI}
\frac{1}{p^n}P_{p,I}(x,x)=\frac{1}{(2\pi)^n}|\mathcal K_I|\sqrt{\det B_x}+\mathcal O(p^{-1}), \quad x\in X,\quad p\to \infty,
\end{equation}
where $\mathcal O$ is uniform on $x\in X$ in $C^k$-topology for any $k$. 

As a consequence, we have
\begin{multline}\label{e:Demailly1}
N_p=\int_XP_{p,I}(x,x)dv_X(x) \\ = \frac{p^n}{(2\pi)^n}|\mathcal K_I|{\rm Vol}_{\mathbf B}(X)+\mathcal O(p^{-1}), \quad p\to \infty.   
\end{multline}

By \eqref{e:def-fp}, \eqref{e:PpI} and \eqref{e:Demailly1}, we infer that
 \[
f_p(x)
={\rm Vol}_{\mathbf B}(X)f_{\mathbf B}(x)+r_p(x),
\]
where
\[
r_p(x)=\mathcal O(p^{-1}), \quad p\to \infty,
 \]
 and $\mathcal O$ is uniform on $x\in X$ in $C^k$-topology for any $k$. 

It follows that 
\begin{multline}
|(f_p(x)-f_p(y))^2-({\rm Vol}_{\mathbf B}(X))^2(f_{\mathbf B}(x)-f_{\mathbf B}(y))^2|\\ \leq C|r_p(x)-r_p(y)|=\mathcal O(p^{-1})d(x,y),\quad x,y\in X, 
\end{multline}  
where $d(x,y)$ is the geodesic distance function defined by the Riemannian metric $g$. 
We infer that
 \[
 \|[P_{p,I},f_p]\|_{2}=\int_X\int_X |P_{p,I}(x,y)|_{h^{L^p}}^2(f_p(x)-f_p(y))^2dv_X(x)dv_X(y)
 \]
 \[
 =({\rm Vol}_{\mathbf B}(X))^2\int_X\int_X |P_{p,I}(x,y)|_{h^{L^p}}^2(f_{\mathbf B}(x)-f_{\mathbf B}(y))^2 dv_X(x)dv_X(y)
 \]
 \[
+\mathcal O(p^{-1}) \int_X\int_X |P_{p,I}(x,y)|_{h^{L^p}}^2d(x,y)dv_X(x)dv_X(y). 
 \]
By \cite{Kor22},Theorem 2], there exist $c>0$ and $C>0$ such that for any $p\in \mathbb N$, $x, y \in X$, we have
\[
|P_{p,I}(x,y)|_{h^{L^p}}\leq Cp^{n}e^{-c\sqrt{p} \,d(x, y)}.
\]
Using this estimate, one can easily show that 
 \[
 \int_X\int_X |P_{p,I}(x,y)|_{h^{L^p}}^2d(x,y)dv_X(x) dv_X(y)=\mathcal O(p^{-1/2}), \quad p\to \infty.
 \]
It follows that 
\begin{multline}\label{e:d2F2}
\|[P_{p,I},f_p]\|_{2}\\
 =({\rm Vol}_{\mathbf B}(X))^2\int_X\int_X |P_{p,I}(x,y)|_{h^{L^p}}^2(f_{\mathbf B}(x)-f_{\mathbf B}(y))^2 dv_X(x)dv_X(y) +\mathcal O(p^{-3/2}). 
\end{multline}  
 
By \cite[Proof of Theorem 1.2]{dpp}, we know that  
\begin{multline}\label{e:VNpft1}
\int_X \int_X |P_{p,I}(x, y)|^2_{h^{L^p}} (f_{\mathbf B}(x)-f_{\mathbf B}(y))^2\, dv_X(x)\, dv_X(y)\\ =\frac{1}{2\pi}\frac{p^{n-1}}{(2\pi)^{n-1}}\int_X|df_{\mathbf B}(x)|_{I}^2\,\Omega_{\mathbf B}(x)+ o(p^{n-1}), \quad p\to \infty. 
\end{multline} 

Combining \eqref{e:VNpft-rhs}, \eqref{e:com-Ppit22}, \eqref{e:com-Ppit}, \eqref{e:d2F2}, \eqref{e:VNpft1}, we get
\begin{multline}\label{e:VNpft}
\int_X \int_X |P_{p,I,t}(x, y)|^2_{h^{L^p}} (f_p(x)-f_p(y))^2d\mu_t(x)d\mu_t(y)\\ =\frac{1}{2\pi}\frac{p^{n-1}}{(2\pi)^{n-1}}\int_X|df_{\mathbf B}(x)|_{I}^2\,\Omega_{\mathbf B}(x)+ o(p^{n-1}), \quad p\to \infty. 
\end{multline}
This is an analog of \cite[Theorem 2.3]{LB24}, which is in its turn an extension of \cite[Theorem 5.8]{Berman18}.
By \eqref{e:d2F1}, \eqref{e:Demailly1}, and \eqref{e:VNpft}, we conclude that
 \begin{equation}\label{e:limd2Fp}
 \lim_{p\to \infty} \frac{d^2F_p}{dt^2}(t)=-\sigma^2, \quad t\in \mathbb R, 
\end{equation}
uniformly on compacts of $\mathbb R$, where $\sigma^2$ is given by \eqref{e:mcarlo}.
 
Now we complete the proof, repeating the arguments of \cite{Berman18,LB24}. First, taking into account \eqref{e:Fp0} and \eqref{e:dFp0}, we can write 
 \[
F_p(t)=\int_0^t\int_0^s \frac{d^2F_p}{dt^2}(u)\,du\, ds,\quad t\in \mathbb R. 
 \]
Using \eqref{e:limd2Fp}, we infer that, for any $t\in \mathbb R$, there exists the limit
 \begin{equation}\label{e:limFp}
 \lim_{p\to \infty}F_p(t)=-\frac{t^2}{2}\sigma^2.
 \end{equation}
We observe that $F_p$ clearly extends to a holomorphic function on $\mathbb C$, and it is uniformly bounded on any compact subset of  $\mathbb C$. Montel's theorem states that a family of locally bounded holomorphic functions is normal, that is, it is a pre-compact subset of the space of continuous functions with respect to compact-open topology. It follows that there exists a subsequence $F_{p_k}$ that converges uniformly to some holomorphic function $F_\infty$ on any compact of $\mathbb C$. By \eqref{e:limFp}, all these limits coincide on $\mathbb R$. From the analytic extension Theorem, we obtain that these limits also coincide on $\mathbb C$, hence we have the uniform convergence $F_p\to F_{\infty}$ on all compacts subsets of $\mathbb C$. If we restrict the previous convergence to the imaginary line $i\mathbb R$, we obtain the convergence of the characteristic function of $\Xi_p$ to the characteristic function of the Gaussian distribution $N(0,\sigma^2)$.


\begin{thebibliography}{00}

\bibitem{B59}
N. S. Bakhvalov. Approximate computation of multiple integrals. (Russian) \textit{Vestnik Moskov. Univ. Ser. Mat. Meh. Astr. Fiz. Him.}  \textbf{1959} (1959), no, 4, 3--18; English translation: On the approximate calculation of multiple integrals, \textit{J. Complexity} \textbf{31} (2015), 502--516.

\bibitem{BH20}
Bardenet, R.; Hardy, A. Monte Carlo with determinantal point processes. \textit{Ann. Appl. Probab.} 30 (2020), no. 1, 368--417.   

\bibitem{Berman14} 
R.J. Berman, Determinantal point processes and fermions on complex manifolds: large deviations and bosonization. \textit{Comm. Math. Phys.} \textbf{327} (2014), 1--47. 

\bibitem{Berman18} 
R.J. Berman, Determinantal point processes and fermions on polarized complex manifolds: bulk universality, In \textit{Algebraic and analytic microlocal analysis,} Springer Proc. Math. Stat., vol. 269, Springer, Cham, 2018, pp. 341--393.

\bibitem{charles20b}
L. Charles, On the spectrum of nondegenerate magnetic Laplacians. \textit{Anal. PDE} \textbf{17} (2024), no. 6, 1907--1952.

\bibitem{GBV19}
G. Gautier, R. Bardenet, M. Valko, On Two Ways to Use Determinantal Point Processes for Monte Carlo Integration. In \textit{Advances in Neural Information Processing Systems (NeurIPS),} volume 32. 2019. Curran Associates, Inc.

 \bibitem{Gu-Uribe}
V. Guillemin and A. Uribe, The Laplace operator on the $n$th tensor power of a line bundle: eigenvalues which are uniformly bounded in $n$. \textit{Asymptotic Anal.} \textbf{1} (1988), 105--113.

\bibitem{HW19}
A. Haimi and A. Wennman, A central limit theorem for fluctuations in Polyanalytic Ginibre ensembles. \textit{Int. Math. Res. Not. IMRN,} \textbf{2019}, no. 5,  1350--1372.

\bibitem{HMBL23}
D. Hawat, G. Mastrilli, R. Bardenet, and R. Lachi\`eze-Rey, Repelled point processes with application to numerical integration, Preprint arXiv:2308.04825. 

\bibitem{Kor18}
Yu. A.  Kordyukov, On asymptotic expansions of generalized Bergman kernels on symplectic manifolds. (Russian); translated from \textit{Algebra i Analiz} \textbf{30} (2018), no. 2, 163--187 \textit{St. Petersburg Math. J.} \textbf{30} (2019), no. 2, 267--283

\bibitem{Kor22}  
Yu. A.  Kordyukov, Semiclassical spectral analysis of the Bochner-Schr\"odinger operator on symplectic manifolds of bounded geometry, \textit{Anal. Math. Phys.} \textbf{12} (2022). no. 1, Paper No 22, 37 pp.

\bibitem{dpp}
Yu. A.  Kordyukov, {Determinantal point processes associated with the Bochner-Schr\"odinger operator}, to appear in \textit{Russ. J. Math. Phys.} \textbf{33} (2026),  preprint  arXiv:2605.13575 (2026). 

\bibitem{Lemoine22}
T. Lemoine. Determinantal point processes associated with Bergman kernels: construction and limit theorems, preprint arXiv:2211.06955 (2022).

\bibitem{LB24}
T. Lemoine and R. Bardenet. Monte Carlo methods on compact complex manifolds using Bergman kernels. preprint arXiv:2405.09203, 2024.

\bibitem{ma-ma:book}
X. Ma and G. Marinescu, \textit{Holomorphic Morse inequalities and Bergman kernels.} Progress in Mathematics, 254. Birkh\"auser Verlag, Basel, 2007. 

\bibitem{ma-ma08} 
X. Ma and G. Marinescu, Generalized Bergman kernels on symplectic manifolds.\textit{Adv. Math.} \textbf{217} (2008), 1756--1815.

\bibitem{N16}
E. Novak. Some results on the complexity of numerical integration. Monte Carlo and quasi-Monte Carlo methods, 161--183, \textit{Springer Proc. Math. Stat., } 163, Springer, [Cham], 2016.

\bibitem{PBD}
V. Petrovic, R. Bardenet, A. Desolneux,
Repulsive Monte Carlo on the sphere for the sliced Wasserstein distance, \textit{Transactions on Machine Learning Research Journal,} In press; Preprint arXiv:2308.04825. 
\end{thebibliography}
\end{document}